\documentstyle[12pt,aasms4]{article}

\newcommand{\SS}{\renewcommand{\baselinestretch}{1}   \tiny \normalsize}

\font\MyScript=eusm10

\begin{document}

\SS
\large

\begin{center}
{\bf  Bayesian Blocks in Two or More Dimensions: \\
Image Segmentation and Cluster Analysis}\\
\vskip 0.5in
{\bf Jeffrey D. Scargle}\\
{\bf NASA Ames Research Center}\\
\end{center}

\vskip 1.75in

\noindent
Contribution to {\bf Workshop on Bayesian Inference and Maximum Entropy Methods in
Science and Engineering (MAXENT 2001)}, held at Johns Hopkins University, 
Baltimore, MD USA on August 4-9, 2001.

\tableofcontents

\section{Summary}

This paper describes an extension, to higher dimensions, of the 
Bayesian Blocks algorithm  for estimating 
signals in noisy time series data (Scargle 1998, 2000). 
The mathematical problem is to find the partition of the data space 
with the maximum posterior probability for
a model consisting of a homogeneous Poisson process 
for each partition element.
For model {\MyScript M$_{n}$}, attributing the data 
within region $n$ of the data space to a Poisson process with a
fixed event\footnote{Throughout we use the terms {\it event} and
{\it(data) point} synonymously.  The event is the occurrence of a 
datum at a given point in the data space.}
 rate $\lambda_{n}$, the global posterior is:
\begin{equation}
P( \mbox{\MyScript M$_{n}$} ) = \Phi(N,V) = { \Gamma( N + 1 ) \Gamma( V - N + 1 ) \over \Gamma( V + 2 ) } = 
{N! (V-N)! \over (V+1)!} \ .
\label{posterior}
\end{equation}
\noindent
Note that $\lambda_{n}$ does not appear, since it has
been marginalized, using a flat, improper prior.
Other priors yield similar formulas.
This expression is valid for 
a data space of any dimension.
It depends on only 
$N$, the number of data points within the region, 
and $V$, the volume of the region. 
No information about the actual locations of the points
enters this expression.

Suppose two such regions, 
described by $N_{1},V_{1}$ and $N_{2},V_{2}$,
are candidates for being merged into one.
From Eq. (\ref{posterior}) construct a {\it Bayes merge factor},
giving the ratio of posteriors for 
the two regions merged and not merged, respectively:
\begin{equation}
P( \mbox{Merge} ) = { \Phi(N_{1}+N_{2},V_{1}+V_{2}) \over
\Phi(N_{1},V_{1})  \Phi(N_{2},V_{2}) } \ .
\label{merge}
\end{equation}
Then  collect data points into {\it blocks} with
this {\it cell coalescence algorithm}:

\begin{itemize}
\item[(0)] Start with the Voronoi tessellation of the data points

\item[(1)] Identify each Voronoi cell as a block 

\item[(2)] Iteratively merge the pair of blocks with the largest 
           {\it merge factor}

\item[(3)] Halt when the maximum merge factor falls below 1

\end{itemize}

\noindent
In many applications it is both required and efficient to
place the restriction that only blocks touching each other
are allowed to merge.
This algorithm partitions the data space into a set of blocks, 
typically much fewer in number than the data points.
Each block has a density equal to the number of data points 
in it divided by its volume.
Then, if desired, high-density blocks adjacent to each other
can be collected into {\it clusters}.

This method allows detection of clusters in high-dimensional 
data spaces, with the following properties:
\begin{itemize}

  \item The number of clusters is determined, not assumed

  \item Clusters can have any shape:

      \begin{itemize}
     \item Avoid the conventional Gaussian assumption 
     \item Shapes can include both concavities and convexities 
     \item Blocks and clusters do not even have to be simply connected
      \end{itemize}

   \item The {\it density profiles} within clusters are estimated, not just the
    locations of the cluster boundaries

  \item Any slowly varying {\it background} is automatically identified

  \item No binning of the raw data is necessary

\end{itemize}

\section{Statistical Challenges of Modern Astronomy}

Certainly in astronomy, and in many other ways, 
this is the {\it Century of Data}.
Several new observational programs 
are severely challenging known techniques for acquiring,
archiving, reducing, analyzing, and interpreting 
astronomical data.
The methods described here have been 
inspired by the need for automated 
extraction of scientific results from 
large, synoptic data sets that will be 
produced by 
projects such as the Sloan Digital Sky Survey,
various other cosmological projects, and -- closer to home --
the exploration of Mars, the Solar System, and planetary
systems in the Sun's neighborhood.

Automated processing already plays a large role
in astronomical data analysis, and will clearly 
become more and more important.
What is not generally agreed on is how far along the path
to the final scientific output automatic processing can be taken.
In my opinion, artificially intelligent processes
will soon become surprisingly practical in this setting.

\section{The Data: Points, Counts, Measurements}

Consider data obtained to estimate
some quantity (as opposed to a discrete attribute).
Such measurements are almost always corrupted by noise, 
blurring, and other instrumental effects.  
The observations may be in a space of one dimension ({\it e.g.}, time series, 
energy spectra), in two dimensions ({\it e.g.} images), or in a space of
higher dimension ({\it e.g.} galaxy redshift/position catalogs).

Three types of measurement data can be distinguished.
The first is {\it point data}, often called {\it event data} in 1D.
That is, it consists of discrete points in the data space under consideration.
The density of points detected in a specific region is 
taken as an estimate of the true density there.
Examples from NASA's Compton $\gamma$-Ray Observatory:
time-tagged photon data,
consisting of lists of photon detection times in units of 2 microseconds;
and sky-image data, consisting of lists of photon positions,
energies and times.

While the usual coordinate representation of such points uses
real numbers, in practice the corresponding infinite accuracy 
and resolution is not achievable.  
The coordinate is usually quantized in some unit, small 
compared to the total range of observation.  
In high energy (x- and $\gamma$-ray) time series data, for example,
the points are the times of detection of individual photons; 
the corresponding quantum is the resolution of the spacecraft clock,
typically in the range of microseconds to milliseconds.

In the second data type 
the entire observation interval (or area, or volume)
is partitioned into pre-specified bins (or pixels, or cells)
and the number of events in each is recorded.
Event data can be converted to this mode,
but in doing so one discards information, reduces the resolution 
of the data, and makes the results 
dependent on the sizes and locations of the bins.
Hence this mode should be used only if event data 
are not available.

A third type of sequential data 
is the measurement of some quantity, at a set of times or points in space.  
This measurement need not be an event count,  
but can be any measurement operation that yields a real number.  
An important difference from the event-counting mode arises 
in the statistical distribution of the 
observational errors: counts in bins typically obey a Poisson distribution,
whereas here the errors may have any distribution, 
often assumed normal (Gaussian).

\section{The Model: Segmentation Yields Structure}
\label{model}

In any analysis involving likelihoods,
a key step is the
choice of a model representing two separate 
aspects of the data: the underlying process, or true signal;
and the noise process corrupting the observations, 
thus hiding the true signal.
Precisely, we must compute the probability, or likelihood, 
that the observed data would be obtained, given the model
and its parameters. 
Although we are studying the underlying signal, 
the likelihood depends also on the data mode,
the sampling process,
and the nature of the noise and other corruption processes.

Data consisting of independent points 
are efficiently described by a very simple model,
the {\it Poisson process}.
In such a process the probability that an event will occur 
in a small element of the data space 
is just a coefficient times the volume of that element.
This coefficient is called the local event rate or {\it Poisson parameter}.
It need not be constant, 
but can vary in an arbitrary way over the data space.
If this variation is random -- a {\it Cox} or 
{\it doubly-stochastic process} -- 
{\bf it is important to distinguish the two random processes.  
One describes the occurrence events at a given location, 
the other how the event rate varies with location.}
If, as usually assumed, these two processes are independent of each other,
even if the event rates at different locations  
are strongly correlated, the occurrence of events 
at these locations are independent of each other 
({\it i.e.}, the joint probability is the product 
of the individual probabilities).

The assumed independence of the data means that the occurrence of one event 
does not affect the probability of any other event.
Hence, a probability referring to separate
subsets of the data space is the product of the
probabilities of the subsets.
A common example of failure of this property -- 
{\it i.e.} dependence -- is ``dead time'' in time series
data: each photon is followed by an interval in which the
occurrence of a second photon is inhibited.

In general, 
the Poisson model is surprisingly appropriate for astronomical processes.
The model seems extremely specialized 
but is in fact remarkably general and serves for almost all applications.  
All we have assumed is that 
the events are independent, and the local event rate 
is specified by an arbitrary density function.
The main limitation, 
need for independence of the events, 
can be relaxed by simply incorporating the 
dependences into the likelihood function.

The model must specify the 
Poisson rate as a function of location within the space,
either parametrically or nonparametrically.
Here we assume one does not have a detailed signal model, 
indicating nonparametric methods.
In particular, consider {\it piecewise constant representations}
of the signal.
This very convenient class of model has the following
properties
\begin{itemize}
\item nonparametric: the number of parameters depends on the number of data 
points\footnote{Rissanen (1989) discusses this seemingly paradoxical definition.
Examples are polynomials, Fourier series, and wavelet expansions.
The idea is that one is representing the structure generically, 
in terms of basis functions, the number of which depends on how
much information is present -- rather than fitting a predefined
shape to the data.}
\item general: capable of representing any reasonable signal
\item simple: underlying rate is constant on finite intervals
\item useful: any physically significant signal property can be computed
  \begin{itemize}
  \item pulse widths, rise times, decay times, amplitudes
  \item background level
  \end{itemize}
\item easy to compute: Bayesian changepoint detection is very effective
\item easily extended to two dimensional and higher data spaces
\item related to classification, cluster detection, and density estimation
\item data adaptive, {\it i.e.} responds to local irregularities
\end{itemize}
\noindent

One can think of this representation as a density
estimator in which blocks take on the role of bins.
The sizes and locations of the bins are not fixed, but are 
determined by the data.
The idea is as follows.
We partition the data space into subsets, or blocks.
The data points within each block are taken 
to be independent and to obey a Poisson distribution
with a constant rate parameter.
Different blocks have different event rates.
In short, we construct a piecewise constant (``blocky'') 
signal approximation.

Our aim is {\it to detect and characterize 
all of the signal structure supported
by the data} -- {\it i.e.} the statistically significant variations.
How best to do this depends on 
the ultimate purpose of the analysis.
Frequently one is interested, 
not in the true signal shape itself, but in 
quantities describing local structures -- {\it e.g.}
pulse widths, amplitudes, {\it etc.}.
Since there are very convenient ways to estimate such parameters
directly from blocks, our seemingly crude representation
may be perfectly adequate.
For such purposes there is no motivation
to impose smoothness, 
although such cosmetic properties are 
important for visualization of model-data relationships.

Note that we don't assert that the true event rate 
changes in a blocky, discontinuous way.  
The physical signal is presumably continuous, 
and we represent it as piecewise constant 
in the spirit of a step-function approximation 
of a smooth curve -- and not with any hope that 
this representation is exact in the limit of small steps
({\it cf.} wavelet theory, and see especially the innovative 
ideas in Donoho 1994a,b).
In the interests of accuracy,
one might consider models of higher complexity,
such as a piecewise {\it linear} representation; 
however, imposing continuity limits the number of free parameters 
(to be about the same as for piecewise constant models), 
so the added accuracy is somewhat illusory.

See (Stoyan, Kendall and Mecke 1995) for an excellent
discussion of point processes in general,
Poisson point processes in particular, and a number of 
ways that real world data can depart from being Poisson.

\section{Algorithms in One Dimension}
\label{algo}

We begin by considering the one dimensional case
(Scargle 1998, 2000).  Time series data are usually considered to
consist of a signal, the character of which is
under investigation, corrupted by 
observational noise\footnote{A point of occasional
confusion is that {\it noise}
in astronomy and physics has two quite distinct meanings:
random observational errors, and random variability intrinsic to
the source or physical system.  The latter, part of the signal, 
is often just what one is studying, whereas the observational
noise is a corruption, to be eliminated as much as possible.}.
In this setting,  segmentation of the data space
into subsets where the signal is taken as constant
is a very practical representation for many signals.
The very different problem of
searching for periodic signals will not be discussed here (see Bretthorst 1988,
more recent publications on his web site, and
my other contribution to this volume).

Three algorithms for implementing this approach
to modeling time series have been described elsewhere
(Scargle 1998, 2000), so I omit details.  
Broadly, the approaches are:
\begin{itemize}
\item {\bf Divide and Conquer}: 
use model comparison to decide whether
the interval should be subdivided; apply iteratively to sub-intervals
\item {\bf Cell Coalescence}: start from an ultra-fine representation
assigning one block to each datum; merge pairs of blocks
based on model comparison
\item {\bf Markov Chain Monte Carlo}: search the block 
parameter space, computing the posterior
(with  all block edge locations and Poisson rate parameters marginalized)
as a function of the number of blocks
\end{itemize}
\noindent
The first two, {\it top-down} and {\it bottom-up}, 
are greedy\footnote{ {\it Greedy algorithms} implement {\it myopic optimization} 
with a ``take what you can now, no regard for the future'' strategy.
On termination of the algorithm 
the resulting local optimum may be a good approximate solution, 
but is not guaranteed to be the global optimum.}
algorithms, iterated until
the answer to the question ``subdivide?'' or ``merge?'' 
is always ``no.''
MCMC convergence is more subtle.

I have been developing the Cell Coalescence method
for astronomical applications, primarily due to its
easy generalization to higher dimensions,
as we shall now discuss.


\section{Structure in Higher Dimensions}

This section describes an iterative algorithm
for structure analysis based on Bayesian methods.
The problem of estimating structure in higher dimensional data spaces
is closely related to {\it regression}, {\it density estimation} (Scott 1992), 
{\it cluster analysis} (Backer 1995), {\it data classification} (Gordon 1999), 
{\it projection pursuit} (Friedman and Tukey 1974), {\it etc.},
depending on the exact nature of the data and the goals of the analysis.
Common to all of these approaches is the same situation 
we have just discussed: the data reflect 
an underlying signal subject to noise and possibly other corruptions;
the signal and the specific noise contributions are unknown, 
but a statistical model of the noise is known or assumed.
The goal is to recover as much information as possible 
about the signal, making good use of any prior information
about the signal and the noise.

Bayesian methods are quite powerful in this arena,
due to the natural way prior information and 
nuisance parameters are handled.
In particular, Bayesian methods provide good solutions
to problems which apparently plague cluster analysis
and related fields 
(Kaufman and Rousseeuw 1990, Backer 1995, Day 1990): 
a bewildering variety of {\it ad hoc}
methods, loss of information in the analysis, 
and assessing {\it post facto} validity of clusters.


\subsection{Voronoi Cells}

The algorithm I am developing for
2D problems uses the same posterior as in
Equations (\ref{posterior})-(\ref{merge}),
since their derivations apply to Poisson
models for point data in any dimension.
However, geometrical considerations for partitioning the data space 
are considerably harder in 2D.
Generalizing the concept of intervals is not very productive, 
due to the great freedom one has in the choice for block shapes:
$\dots$ squares? rectangles? circles? other?
What {\it is} productive is to follow the cell coalescence idea
(above, \S \ref{algo})
and construct small, 
elementary blocks from which larger blocks can be composed: 
{\it assemble the macroscopic out of the microscopic}.

In the 1D case, the intervals spanned by the midpoints 
between successive pairs of adjacent data points are
the obvious choice for the cells.
The mid-point based construction has the following properties:
\begin{enumerate}
\item[(1)] Cell $i$ contains all points closer to datum $i$ than to any other datum
\item[(2)] Forms a {\it partition}: Every point in the space lies in one and only one cell 
\item[(3)] Attached to each cell is a number, $A_{i}$, namely its area
\item[(4)] ${1 \over A_{i} }$ is a measure of the density of data points in the vicinity of datum $i$
\end{enumerate}
\noindent 
Item (1), a constructive definition valid in any dimension,
is the space's {\it Voronoi tessellation} 
({\it e.g.}, Stoyan, Kendall and Mecke 1995) 
determined by the data points.  
This definition implies (2), {\it i.e.} Voronoi cells 
{\it partition} the space.  
Where points are close together the cell areas are small, 
and {\it vice versa}.  
Specifically, the reciprocal of the cell area is a convenient measure 
of local point density.
How useful this representation is!  With almost no work
we already have a detailed, if choppy, density estimation.

Now three definitions.  A {\it cell} is one of the Voronoi cells, 
which are in a one-to-one correspondence with the data points.  
During the iteration we insure that these cells maintain their identity, 
so there is no problem thinking of cells and data points interchangeably.

A {\it block} is a set of one or more merged cells,
including both the initial one-cell blocks and the multi-cell
blocks that form during the iteration.
It is common to restrict mergers to, {\it e.g.}, only blocks that are 
touching\footnote{That two blocks are touching can be readily determined by 
finding a Voronoi cell in one block that has a vertex in common 
with at least one cell in the other block.
The MatLab (\copyright \ the Math Works, Inc.) routine computing 
Voronoi tessellations in arbitrary
dimensions returns vertex information in a form
very convenient for all such computations.};
hence blocks consist of adjacent cells.

One last definition: At the end of the iteration, we may find that some sets of adjacent 
blocks stand out from the general background, 
say by having event rates 
substantially above background.  
We call such groups {\it clusters}.  
In the case of galaxies, this corresponds precisely
to the notion of {\it galaxy clusters}.
In the case of photons, this corresponds 
to the notion of a {\it source}.
Note that we are interested in density profiles within clusters;
for the most part traditional classification and cluster analysis 
identifies boundaries only.

\subsection{The Algorithm: Cell Coalescence}

We are now ready to outline the algorithm.
The idea is simple: 
find the partition with the largest total posterior probability
\begin{equation}
P_{total} \equiv \Pi_{n=1}^{N} \Phi(N_{n},V_{n}) \ ,
\label{posterior_sum}
\end{equation}
\noindent
obtained from Eq. (\ref{posterior}) 
applied to each element of the partition, and using the
product rule since the elements are assumed independent.
In principle we want to maximize this posterior over the set 
\begin{equation}
S_{\infty} \equiv \{ \mbox{\MyScript{P}} \ | \  
 \mbox{ D is partitioned by } \mbox{\MyScript{P}} \}
\label{all}
\end{equation}
\noindent
where $D$ is the entire data space.
This task would be awkwardly complicated and impractically large,
so instead we restrict attention to the set
\begin{equation}
S_{0} \equiv \{ \mbox{\MyScript{P}} \ | \   
\mbox{ D is partitioned by \mbox{\MyScript{P}} into blocks of Voronoi cells } \} \ .
\label{some}
\end{equation}
\noindent
There are many ways in which 
these two sets are different  -- (\ref{all}) is highly infinite,
(\ref{some}) is quite finite.
But for representing the structural information 
content of the data points $S_{\infty}$ and $S_{0}$ are effectively the same.  
This is sometimes expressed 
by saying that the Voronoi tessellation contains 
all the proximity information in the data,
a well known fact in the computer graphics industry.

I have coded and experimented with a very 
simple search through $S_{0}$, namely
a greedy iteration that at each step merges the pair of blocks with the
largest value of the merge posterior in Eq. (\ref{merge}):
\begin{itemize}
\item[(0)] Initialize: one block ($\equiv$ a Voronoi cell) per data point
\item[(1)] Compute change in $P(\mbox{\MyScript{P}})$ for merging 
   each pair of adjacent blocks
\item[(2)] If largest change is positive, merge\footnote{To 
merge cells $n$ and $m$ simply change cell $n$ according to
$N_{n} = N_{n} + N_{m}, V_{n} = V_{n} + V_{m}$ and delete cell $m$.}
 that pair and Go To (1)
\item[(3)] Else (no pair should be merged) stop.
\end{itemize}
\noindent
This greedy algorithm can be criticized because it is 
not guaranteed to yield the global maximum.
To assess this problem, I am exploring alternative procedures 
which explore more of the parameter space -- including randomized 
algorithms (Motwani and Raghavan 1995),
algorithms using more general transformations than just block merges,
and combinatorial optimization methods based on dynamic programming
(Hubert, Arabie, and Meulman 2001).

\subsection{Adopting a Model}

There are several ways to determine the model to be adopted. 
Most simply one can use the model in place 
when the stopping criterion is triggered.  
Or, one can select the largest $P(\mbox{\MyScript{P}})$ from 
the history of the iteration.  
If the transformations are such that $P(\mbox{\MyScript{P}})$ is monotonically non-decreasing, 
these two are obviously the same.  An alternative is to 
average all of the models, weighting each one according to its
posterior probability -- {\it model averaging}
(see Hoeting, Madigan, Raftery, and Volinsky 1999,
and for some software \verb+http://www.research.att.com/~volinsky/bma.html+).
This approach has the advantage of 
producing a model that is more in accord with the 
maximum a posteriori (MAP) principle, but the disadvantage that the
averaged model is not in the initially defined space of models -- namely, 
piecewise constant Poissonian blocks.
Depending on what information one wants from the models, 
this disadvantage may or may not be important.

\subsection{Related Work}

Various authors 
(Ebeling and Wiedenmann 1993, Scargle 2000,
Ramella, Boschin, Fadda and Nonino 2001) 
have considered algorithms that start from the Voronoi tessellation 
and iteratively merge the Voronoi cells into blocks based on
some kind of statistical criterion.
Essentially the same problem has been discussed in the literature of automatic data 
classification (Gordon 1999, Backer 1995), where the concept of merging cells is 
replaced by a more generalized set of allowed transformations of the partition.

Our ultimate goal is the same as to that of classification, 
the ``detection of important relationships and 
structure within data'' (Gordon 1999, \S 1.2) by finding ``a partition in 
which objects are similar to the other objects belonging to their class and 
dissimilar to objects belonging to different classes'' (Gordon 1999, \S 3.1).  
However, we generalize the concept of similarity 
(most often taken to be simply closeness of the points)
to mean that the distributions of the observed 
points are adequately described by an appropriate statistical model.
Other classification criteria, such as 
homogeneity, lack of heterogeneity, 
or isolation or separation from the rest of the data, 
can be subsumed into this more general concept.

\section{Information Theoretic Inference}

Several foundational issues are pertinent to the above analysis,
beginning with the Bayesian definition of 
{\it probability} as a measure
of the degree of one's belief in the truth of
the hypothesis (Jaynes 1997).  It would be hard to identify
two words more fraught with philosophical difficulties 
than {\it belief} and {\it truth}.

Science never achieves absolute truth of any hypothesis.
Instead it seeks to establish approximate but useful descriptions,
models, predictions of future events and measurements, and so on.
Belief -- by others or of one's own -- is 
hard to define operationally.
The pragmatic point here is that 
{\bf we believe that our model is not true}.
A discontinuous model can not be a true 
representation of a continuous signal.
The prior for our model is zero.
This leads directly to oblivion.
There are ways to deal with this problem
within the Bayesian framework, 
but uneasiness about this point remains.

The vexing concept of the {\it noninformative prior}
is particularly relevant for the construction of
general purpose, turn-key algorithms.  
Although it is possible to fashion algorithms
so that one can plug in the prior that expresses
one's knowledge in any given situation, it is still
legitimate to seek an algorithm that makes no prior
assumptions, perhaps at the price of some loss of
sensitivity, efficiency, or the like.

There are other aspects of the Bayesian
data analysis formalism that are muddled 
by the belief-in-truth interpretation.
The founders of Bayesian data analysis
made a great step forward in recognizing that
{\it deductive} logic is inappropriate to scientific
thought and research.
However, in substituting 
{\it inductive} logic 
they missed that 
the problem is not just the kind of logic,
but that logic itself is inappropriate
to the practice of scientific data analysis.

The notion that science deals with information
about the Universe, more than with truths about it, leads
to a data analysis approach that substitutes information theoretic
quantities for the subjective Bayesian probabilities.
In particular {\it information contained in data and models}
replaces the usual notion of {\it belief in the truth of hypotheses}.
The resulting framework 
avoids difficulties associated with noninformative priors, 
and I believe is more in tune with how data analysis
actually proceeds -- uncovering and processing information
about the Universe.

We use models to extract from data information about the world.
There is no learning of absolute truths about the world.
We proceed by comparing model predictions
with observed data.
Evaluating the joint information contained
in model and data is a straightforward application
of information theory.  
Elsewhere in this volume,
related information theoretic ideas are
discussed in the papers
``Inductive Logic'' by Kevin Knuth, and 
``Information Theoretic Approach to Bayesian Inference''
by Jeffrey Jewell.

We equate two expressions
for data-model information -- following the 
standard derivation of Bayes Theorem 
by equating two expressions for the joint data-model probability.
Start by defining the
{\it mutual information} of two random variables 
as the relative entropy between their actual joint distribution 
and what the joint distribution would have been if the
variables were independent -- {\it i.e.} the product of
the individual distributions:

\noindent
{\bf Definition: Mutual Information}\\
\noindent
Given random variables $X$ and $Y$,
with probability distributions $p(x)$ and $p(y)$,
and joint probability distribution $p(x,y)$,
the mutual information of $X$ and $Y$ is
\begin{equation}
I(X;Y) = \sum_{x} \sum_{y} p(x,y) \  log {p(x,y) \over p(x) p(y) }
\end{equation}
\noindent
\noindent
that is, the expectation of $log {p(x,y) \over p(x) p(y) }$
with respect to the full joint distribution $p(x,y)$.
The sums in $x$ and $y$ are over the appropriate 
discrete event spaces, and become integrals for continuous
variables.

Cover and Thomas (1991) prove that
\begin{equation}
I(X;Y) = H(X) - H(X|Y)
\end{equation}
\noindent
where $H(X)$ is the unconditional entropy of $X$,
and $H(X|Y)$ is the entropy of $X$ conditional on $Y$.
Noting that $I(X;Y) = I(Y;X)$, 
using the more mnemonic variable names $D$ for data and $M$ for model,
and replacing entropy, $H$, with the 
negative of information, $-I$, we have
\begin{equation}
I(M|D) = I(D|M) + I(M) - I(D) \ ,
\end{equation}
\noindent
rather like Bayes Theorem with the 
identification $I \leftrightarrow log(P)$.
I am experimenting with inference procedures
for cluster analysis based on this relationship.

\vskip .25in

I am indebted to Tom Loredo, Alanna Connors, Larry Bretthorst,
Peter Sturrock, Jay Norris, Jerry Bonnell, Ayman Farahat,
and other colleagues
for comments, suggestions, and encouragement,
and to Michele Bierbaum and Bob Fry for their kind
hospitality at MAXENT 2001.

{99}

\end{document}